\documentclass[a4paper]{amsart}  

\usepackage{amssymb} \usepackage{amscd} \usepackage{times}

\def\zbb{\mathbb{Z}}  
  
  \def\phi{\varphi}
 \def\p1{{\mathbb{P}^1_\zbb}}

\begin{document}

\title{ Sup-Inf inequality On Manifold of Dimension 3}

\author{Samy Skander Bahoura}

\address{Department of Mathematics, Patras University, 26500 Patras , Greece }
              
\email{samybahoura@yahoo.fr, bahoura@ccr.jussieu.fr} 

\date{}

\maketitle

\begin{abstract}

We give an estimate of type $ \sup \times \inf $ on Riemannian manifold of Dimension 3 for prescribed curvature equation.

\end{abstract}

\begin{center} INTRODUCTION AND RESULTS.

\end{center}

\bigskip

We are on Riemannian manifold $ (M,g) $ of dimension 3 not necessary compact. In this paper we denote $ \Delta = -\nabla^j(\nabla_j) $ the geometric laplacian.

\smallskip

The scalar curvature equation is:

$$ 8\Delta u+R_gu=V u^{N-1},\,\, u >0 . \qquad (E) $$

Where $ R_g $ is the scalar curvature and $ V $ is a function (prescribed scalar curvature).

\bigskip

We consider three positive real number $ a,b A $ and we suppose $ V $ lipschitzian:

 $$ 0 < a \leq V(x) \leq b < + \infty \,\, {\rm and} \,\, ||\nabla V||_{L^{\infty}(M)} \leq A. \qquad (C) $$

\bigskip

The equation $ (E) $ was studied lot of when $ M =\Omega \subset {\mathbb R}^n $ or $ M={\mathbb S}_n $ see for example, [B], [CL1], [L1]. In this case we have a $ \sup \times \inf $ inequality.

\smallskip

The corresponding equation in two dimensions on open set $ \Omega $ of $ {\mathbb R}^2 $, is:

$$ \Delta u=Ve^u, \qquad (E') $$

The equation $ (E') $ was studed lot of and we can find very important result about a priori estimates in [BM], [BLS], [CL2], [L2], and [S].

\smallskip

In the case $ V\equiv 1 $ and $ M $ compact, the equation $ (E) $ is Yamabe equation. It was studed lot of, T.Aubin and R.Schoen have proved the existence of solution in this cas, see for example [Au] and [L-P].

\bigskip

When $ M $ is a compact Riemannian manifold, it exist some compactness result for equation  $ (E) $ see [L-Zh]. Li and Zhu [L-Zh], proved that the energy is bounded and if we suppose $ M $ not diffeormorfic to the three sphere, the solutions are uniformly bounded. To have this result they use the positive mass theorem.

\bigskip

Now, if we suppose $ M $ Riemannian manifold (not necessarily compact) and $ V\equiv 1 $, Li and Zhang [L-Z] proved that the product $ \sup \times \inf $ is bounded.

\bigskip

Here, we give an equality of type $ \sup \times \inf $ for the equation $ (E) $ with general conditions $ (C) $. We have:

\bigskip

{\it Theorem. For all compact set $ K $ of $ M $ and all positive numbers a,b,A, it exists a positive constant c, which depends only on, K,a,b,A,M,g such that:

$$ (\sup_K  u )^{1/3}  \times \inf_M u \leq c, $$

for all $ u $ solution of $ (E) $ with conditions $ (C) $.}

\bigskip

Note that in our work, we have not assumption on energy or boundary condition if we suppose the manifold $ M $ with boundary.

\bigskip

Next, in the proof of the previous theorem, we can replace the scalar curvature by any smooth function $ f $, but here we do the proof with $ R_g $ the scalar curvature.

\newpage

\underbar {\bf Proof of the theorem.}

\bigskip

{\underbar {\bf Part I: The metric and the laplacian in polar coordinates.}}

\bigskip

Let $ (M,g) $ a Riemannian manifold. We note $ g_{x,ij} $ the local expression of the metric $ g $  in the exponential map centred in $ x $.

We are concerning by the polar coordinates expression of the metric. By using Gauss lemma, we can write:

$$ g=ds^2=dt^2+g_{ij}^k(r,\theta)d\theta^id\theta^j=dt^2+r^2{\tilde g}_{ij}^k(r,\theta)d\theta^id\theta^j=g_{x,ij}dx^idx^j, $$

in a polar chart with origin $ x $", $ ]0,\epsilon_0[\times U^k $, with $ ( U^k, \psi) $ a chart of $ {\mathbb S}_{n-1} $. We can write the element volume:

$$ dV_g=r^{n-1}\sqrt {|{\tilde g}^k|}dr d \theta^1 \ldots d \theta^{n-1} = \sqrt {[det(g_{x,ij})]}dx^1 \ldots dx^n, $$

then,

$$ dV_g=r^{n-1} \sqrt { [det(g_{x,ij})]}[\exp_x(r\theta)]\alpha^k(\theta)dr d\theta^1 \ldots d \theta^{n-1} , $$

where, $ \alpha^k $ is such that, $ d\sigma_{{\mathbb S}_{n-1}}=\alpha^k(\theta) d \theta^1 \ldots d \theta^{n-1} . $ (Riemannian volume element of the la sphere in the chart $ (U^k,\psi) $ ).

\bigskip

Then,

$$ \sqrt { |{\tilde g}^k|}=\alpha^k(\theta) \sqrt {[det(g_{x,ij})]}, $$

Clearly, we have the following proposition:

\bigskip

\underbar {\bf Proposition 1:} Let $ x_0 \in M $, there exist $ \epsilon_1>0 $ and if we reduce $ U^k $, we have:

$$ |\partial_r{\tilde g}_{ij}^k(x,r,\theta)|+|\partial_r\partial_{\theta^m}{\tilde g}_{ij}^k(x,r, \theta)| \leq C r,\,\, \forall \,\, x\in B(x_0,\epsilon_1) \,\, \forall \,\, r\in [0,\epsilon_1], \,\, \forall \,\, \theta \in U^k.$$

and,

$$ |\partial_r|{\tilde g}^k|(x,r,\theta)|+\partial_r \partial_{\theta^m} |{\tilde g}^k|(x,r,\theta)\leq C r,\,\, \forall \,\, x\in B(x_0,\epsilon_1) \,\, \forall \,\, r\in [0,\epsilon_1], \,\, \forall \,\, \theta \in U^k. $$

\underbar {\bf Remark:} 

\bigskip

$ \partial_r [ \log \sqrt { |{\tilde g}^k |}] $ is a local function of $ \theta $, and the restriction of the global function on the sphere $ {\mathbb S}_{n-1} $, $ \partial_r [\log \sqrt { det(g_{x, ij})}] $. We will note, $ J(x,r,\theta)=\sqrt { det(g_{x, ij})} $.

\bigskip

Let's write the laplacian in $ [0,\epsilon_1]\times U^k $,

$$ -\Delta = \partial_{rr}+\dfrac{n-1}{r}\partial_r+ \partial_r [\log \sqrt { |{\tilde g^k|}] }\partial_r+\dfrac{1}{r^2 \sqrt {|{\tilde g}^k|}}\partial_{\theta^i}({\tilde g}^{\theta^i \theta^j}\sqrt { |{\tilde g}^k|}\partial_{\theta^j}) . $$

We have,

$$ -\Delta = \partial_{rr}+\dfrac{n-1}{r}\partial_r+ \partial_r \log J(x,r,\theta)\partial_r+ \dfrac{1}{r^2 \sqrt {|{\tilde g}^k|}}\partial_{\theta^i}({\tilde g}^{\theta^i \theta^j}\sqrt { |{\tilde g}^k|}\partial_{\theta^j}) . $$

We write the laplacian ( radial and angular decomposition),

$$ -\Delta = \partial_{rr}+\dfrac{n-1}{r} \partial_r+\partial_r [\log J(x,r,\theta)] \partial_r-\Delta_{{S}_r(x)}, $$

where $ \Delta_{ S_r(x)} $ is the laplacian on the sphere $ {S}_r(x) $. 

\bigskip

We set $ L_{\theta}(x,r)(...)=r^2\Delta_{ S_r(x)}(...)[\exp_x(r\theta)] $, clearly, this operator is a laplacian on $ {\mathbb S}_{n-1} $ for particular metric. We write,

$$ L_{\theta}(x,r)=\Delta_{g_{x,r, {}_{{\mathbb S}_{n-1}}}}, $$

and,

$$ \Delta = \partial_{rr}+\dfrac{n-1}{r} \partial_r+\partial_r [ J(x,r,\theta)] \partial_r - \dfrac{1}{r^2} L_{\theta}(x,r) . $$

If, $ u $ is function on $ M $, then, $ \bar u(r,\theta)=uo\exp_x(r\theta) $ is the corresponding function in polar coordinates centred in $ x $. We have,

$$ -\Delta u =\partial_{rr} \bar u+\dfrac{n-1}{r} \partial_r \bar u+\partial_r [ J(x,r,\theta)] \partial_r \bar u-\dfrac{1}{r^2}L_{\theta}(x,r)\bar u . $$

{\underbar { \bf Part II: "Blow-up" and "Moving-plane" methods }} 

\bigskip

\underbar {\bf The "blow-up" technic}

\bigskip

Let, $ (u_i)_i $ a sequence of functions on $ M $ such that,

$$ \Delta u_i + R_gu_i =V_i{u_i}^5, \,\, u_i>0, \qquad (E')  $$

We argue by contradiction and we suppose that $ \sup^{1/3} \times \inf $ is not bounded.

\smallskip

We assume that:

\bigskip

$ \forall \,\, c,R >0 \,\, \exists \,\, u_{c,R} $ solution of $ (E') $ such that:

$$ R [\sup_{B(x_0,R)} u_{c,R}]^{1/3} \times \inf_M u_{c,R} \geq c, \qquad (H) $$

\bigskip

\underbar {\bf Proposition 2:} 

\bigskip

There exist a sequence of points $ (y_i)_i $, $ y_i \to x_0 $ and two sequences of positive real number $ (l_i)_i, (L_i)_i $, $ l_i \to 0 $, $ L_i \to +\infty $, such that if we consider $ v_i(y)=\dfrac{u_io\exp_{y_i}(y)}{u_i(y_i)} $, we have:

$$ 0 < v_i(y) \leq  \beta_i \leq 2^{1/2}, \,\, \beta_i \to 1. $$

$$  v_i(y)  \to \left ( \dfrac{1}{1+{|y|^2}} \right )^{1/2}, \,\, {\rm uniformly \,\, on\,\, every \,\, compact \,\, set \,\, of } \,\, {\mathbb R}^3 . $$

$$ l_i^{1/2} [u_i(y_i)]^{1/3} \times \inf_M u_i \to +\infty $$

\underbar {\bf Proof:}

We use the hypothesis $ (H) $, we can take two sequences $ R_i>0, R_i \to 0 $ and $ c_i \to +\infty $, such that,

$$ R_i [\sup_{B(x_0,R_i)}]^{1/3} u_i \times \inf_M u_i \geq c_i \to +\infty, $$

Let, $ x_i \in  { B(x_0,R_i)} $, such that $ \sup_{B(x_0,R_i)} u_i=u_i(x_i) $ and $ s_i(x)=[R_i-d(x,x_i)]^{1/2} u_i(x), x\in B(x_i, R_i) $. Then, $ x_i \to x_0 $.

\bigskip

We have, 

$$ \max_{B(x_i,R_i)} s_i(x)=s_i(y_i) \geq s_i(x_i)={R_i}^{1/2} u_i(x_i)\geq \sqrt {c_i}  \to + \infty. $$ 

\bigskip

Set :

$$ l_i=R_i-d(y_i,x_i),\,\, \bar u_i(y)= u_i o \exp_{y_i}(y),\,\,  v_i(z)=\dfrac{u_i o \exp_{y_i}\left ( z/[u_i(y_i)]^{2} \right ) } {u_i(y_i)}. $$

Clearly, $ y_i \to x_0 $. We obtain:

$$ L_i= \dfrac{l_i}{(c_i)^{1/2}} [u_i(y_i)]^{2}=\dfrac{[s_i(y_i)]^{2}}{c_i^{1/2}}\geq \dfrac{c_i^{1}}{c_i^{1/2}}=c_i^{1/2}\to +\infty. $$

\bigskip

If $ |z|\leq L_i $, then $ y=\exp_{y_i}[z/ [u_i(y_i)]^{2}] \in B(y_i,\delta_i l_i) $ with $ \delta_i=\dfrac{1}{(c_i)^{1/2}} $ and $ d(y,y_i) < R_i-d(y_i,x_i) $, thus, $ d(y, x_i) < R_i $ and, $ s_i(y)\leq s_i(y_i) $, we can write,

$$ u_i(y) [R_i-d(y,y_i)]^{1/2} \leq u_i(y_i) (l_i)^{1/2}. $$

But, $ d(y,y_i) \leq \delta_i l_i $, $ R_i >l_i$ and $ R_i-d(y, y_i) \geq R_i-\delta_i l_i>l_i-\delta_i l_i=l_i(1-\delta_i) $, we obtain,

$$ 0 < v_i(z)=\dfrac{u_i(y)}{u_i(y_i)} \leq \left [ \dfrac{l_i}{l_i(1-\delta_i)} \right ]^{1/2}\leq 2^{1/2} . $$

We set, $ \beta_i=\left ( \dfrac{1}{1-\delta_i} \right )^{1/2} $, clearly $ \beta_i \to 1 $.

\bigskip

The function $ v_i $ is solution of:

$$ -g^{jk}[\exp_{y_i}(y)]\partial_{jk} v_i-\partial_k \left [ g^{jk}\sqrt { |g| } \right ][\exp_{y_i}(y)]\partial_j v_i+ \dfrac{ R_g[\exp_{y_i}(y)]}{[u_i(y_i)]^{4}} v_i=V_i{v_i}^5, $$

By elleptic estimates and Ascoli, Ladyzenskaya theorems, $ ( v_i)_i $ converge uniformely on each compact to the function $ v $ solution on $ {\mathbb R}^3 $ of, 

$$ \Delta v=V(x_0)v^5, \,\, v(0)=1,\,\, 0 \leq v\leq 1\leq 2^{1/2}, $$

Without loss of generality, we can suppose $ V(x_0)=3 $.

\bigskip

By using maximum principle, we have $ v>0 $ on $ {\mathbb R}^3 $, the result of Caffarelli-Gidas-Spruck ( see [C-G-S]) give, $ v(y)=\left ( \dfrac{1}{1+{|y|^2}} \right )^{1/2} $. We have the same propreties for $ v_i $ in the previous paper [B2].

\bigskip

\underbar {\bf Polar coordinates and "moving-plane" method}

\bigskip

Let, 

$$ w_i(t,\theta)=e^{1/2}\bar u_i(e^t,\theta) = e^{t/2}u_io\exp_{y_i}(e^t\theta), \,\, {\rm et} \,\, a(y_i,t,\theta)=\log J(y_i,e^t,\theta). $$ 

\underbar {\bf Lemma 1:}

\bigskip

The function $ w_i $ is solution of:

$$  -\partial_{tt} w_i-\partial_t a \partial_t w_i-L_{\theta}(y_i,e^t)+c w_i=V_i w_i^5, $$
 
with,

 $$ c = c(y_i,t,\theta)=\left ( \dfrac{1}{2} \right )^2+ \dfrac{1}{2} \partial_t a  -\lambda e^{2t}, $$

\underbar {\bf Proof:}

\smallskip

We write:

$$ \partial_t w_i=e^{3t/2}\partial_r \bar u_i+\dfrac{1}{2} w_i,\,\, \partial_{tt} w_i=e^{5t/2} \left [\partial_{rr} \bar u_i+\dfrac{2}{e^t}\partial_r \bar u_i \right ]+\left ( \dfrac{1}{2} \right )^2 w_i. $$

$$ \partial_t a =e^t\partial_r \log J(y_i,e^t,\theta), \partial_t a \partial_t w_i=e^{5t/2}\left [ \partial_r \log J\partial_r \bar u_i \right ]+\dfrac{1}{2} \partial_t a w_i.$$

the lemma  is proved.

\bigskip

Now we have, $ \partial_t a=\dfrac{ \partial_t b_1}{b_1} $, $ b_1(y_i,t,\theta)=J(y_i,e^t,\theta)>0 $,

\bigskip

We can write,

$$ -\dfrac{1}{\sqrt {b_1}}\partial_{tt} (\sqrt { b_1} w_i)-L_{\theta}(y_i,e^t)w_i+[c(t)+ b_1^{-1/2} b_2(t,\theta)]w_i=V_i{w_i}^{N-1}, $$

where, $ b_2(t,\theta)=\partial_{tt} (\sqrt {b_1})=\dfrac{1}{2 \sqrt { b_1}}\partial_{tt}b_1-\dfrac{1}{4(b_1)^{3/2}}(\partial_t b_1)^2 .$

\bigskip

Let,

$$ \tilde w_i=\sqrt {b_1} w_i. $$

\underbar {\bf Lemma 2:}

\smallskip

The function $ \tilde w_i $ is solution of:

$$ -\partial_{tt} \tilde w_i+\Delta_{g_{y_i, e^t, {}_{{\mathbb S}_{2}}}} (\tilde w_i)+2\nabla_{\theta}(\tilde  w_i) .\nabla_{\theta} \log (\sqrt {b_1})+(c+b_1^{-1/2} b_2-c_2) \tilde w_i= $$

$$ = V_i\left (\dfrac{1}{b_1} \right )^{2} {\tilde w_i}^{5}, $$

where, $ c_2 =[\dfrac{1}{\sqrt {b_1}} \Delta_{g_{y_i, e^t, {}_{{\mathbb S}_{n-1}}}}(\sqrt{b_1}) + |\nabla_{\theta} \log (\sqrt {b_1})|^2] . $

\bigskip

\underbar {\bf Proof:}

\bigskip

We have: 

$$ -\partial_{tt} \tilde w_i-\sqrt {b_1} \Delta_{g_{y_i, e^t, {}_{{\mathbb S}_{2}}}} w_i+(c+b_2) \tilde w_i= V_i\left (\dfrac{1}{b_1} \right )^{2} {\tilde w_i}^5, $$

But,

$$ \Delta_{g_{y_i, e^t, {}_{{\mathbb S}_{2}}}} (\sqrt {b_1} w_i)=\sqrt {b_1} \Delta_{g_{y_i, e^t, {}_{{\mathbb S}_{2}}}} w_i-2 \nabla_{\theta} w_i .\nabla_{\theta} \sqrt {b_1}+ w_i \Delta_{g_{y_i, e^t, {}_{{\mathbb S}_{2}}}}(\sqrt {b_1}), $$

and,

$$ \nabla_{\theta} (\sqrt {b_1} w_i)=w_i \nabla_{\theta} \sqrt {b_1}+ \sqrt {b_1} \nabla_{\theta} w_i, $$

we deduce,

$$  \sqrt {b_1} \Delta_{g_{y_i, e^t, {}_{{\mathbb S}_{2}}}} w_i= \Delta_{g_{y_i, e^t, {}_{{\mathbb S}_{2}}}} (\tilde w_i)+2\nabla_{\theta}(\tilde  w_i) .\nabla_{\theta} \log (\sqrt {b_1})-c_2 \tilde w_i, $$

with $ c_2=[\dfrac{1}{\sqrt {b_1}} \Delta_{g_{y_i, e^t, {}_{{\mathbb S}_{2}}}}(\sqrt{b_1}) + |\nabla_{\theta} \log (\sqrt {b_1})|^2] . $ The lemma is proved.

\bigskip

\underbar {\bf The "moving-plane" method:}

\bigskip

Let $ \xi_i $ a real number,  and suppose $ \xi_i \leq t $, we set $ t^{\xi_i}=2\xi_i-t $ and $ \tilde w_i^{\xi_i}(t,\theta)=\tilde w_i(t^{\xi_i},\theta) $.

\bigskip

We have, 

$$ -\partial_{tt} \tilde w_i^{\xi_i}+\Delta_{g_{y_i, e^{t^{\xi_i}} {}_{{\mathbb S}_{2}}}} (\tilde w_i)+2\nabla_{\theta}(\tilde  w_i^{\xi_i}) .\nabla_{\theta} \log (\sqrt {b_1}) \tilde w_i^{\xi_i}+[c(t^{\xi_i})+b_1^{-1/2}(t^{\xi_i},.)b_2(t^{\xi_i})-c_2^{\xi_i}] \tilde w_i^{\xi_i}= $$

$$ =V_i^{\xi_i}\left (\dfrac{1}{b_1^{\xi_i}} \right )^{2} { ({\tilde w_i}^{\xi_i}) }^5. $$

By using the same arguments than in [B], we have:

\bigskip

\underbar {\bf Proposition 3:}

\smallskip

We have:

$$ 1)\,\,\, \tilde w_i(\lambda_i,\theta)-\tilde w_i(\lambda_i+4,\theta) \geq \tilde k>0, \,\, \forall \,\, \theta \in {\mathbb S}_{2}. $$

For all $ \beta >0 $, there exist $ c_{\beta} >0 $ such than:

$$ 2) \,\,\, \dfrac{1}{c_{\beta}} e^{t/2}\leq \tilde w_i(\lambda_i+t,\theta) \leq c_{\beta}e^{t/2}, \,\, \forall \,\, t\leq \beta, \,\, \forall \,\, \theta \in {\mathbb S}_{2}. $$

We set,

$$ \bar Z_i=-\partial_{tt} (...)+\Delta_{g_{y_i, e^t, {}_{{\mathbb S}_{2}}}} (...)+2\nabla_{\theta}(...) .\nabla_{\theta} \log (\sqrt {b_1})+(c+b_1^{-1/2} b_2-c_2)(...) $$

{\bf Remark:} In the opertor $ \bar Z_i $, by using the proposition 3, the coeficient $ c+b_1^{-1/2}b_2-c_2 $ verify:

$$ c+b_1^{-1/2}b_2-c_2 \geq k'>0,\,\, {\rm pour }\,\, t<<0, $$

it is fundamental if we want to apply Hopf maximum principle.

\bigskip

\underbar {\bf Goal:}

\bigskip

Like in [B], we have elliptic second order operator, here it's $ \bar Z_i $, the goal is to use the "moving-plane" method to have a contradiction. For this, we must have:

$$ \bar Z_i(\tilde w_i^{\xi_i}-\tilde w_i) \leq 0, \,\, {\rm if} \,\, \tilde w_i^{\xi_i}-\tilde w_i \leq 0. $$

We write:

$$ \bar Z_i(\tilde w_i^{\xi_i}-\tilde w_i)= (\Delta_{g_{y_i, e^{t^{\xi_i}}, {}_{{\mathbb S}_{2}}}}-\Delta_{g_{y_i, e^{t}, {}_{{\mathbb S}_{2}}}}) (\tilde w_i^{\xi_i})+ $$

$$ +2(\nabla_{\theta, e^{t^{\xi_i}}}-\nabla_{\theta, e^t})(w_i^{\xi_i}) .\nabla_{\theta, e^{t^{\xi_i}}} \log (\sqrt {b_1^{\xi_i}})+ 2\nabla_{\theta,e^t}(\tilde w_i^{\xi_i}).\nabla_{\theta, e^{t^{\xi_i}}}[\log (\sqrt {b_1^{\xi_i}})-\log \sqrt {b_1}]+ $$ 

$$ +2\nabla_{\theta,e^t} w_i^{\xi_i}.(\nabla_{\theta,e^{t^{\xi_i}}}-\nabla_{\theta,e^t})\log \sqrt {b_1}- [(c+b_1^{-1/2} b_2-c_2)^{\xi_i}-(c+b_1^{-1/2}b_2-c_2)]\tilde w_i^{\xi_i} + $$

$$ + V_i^{\xi_i}\left ( \dfrac{1}{b_1^{\xi_i}} \right )^{2} ({\tilde w_i}^{\xi_i})^5-V_i\left ( \dfrac{1}{b_1} \right )^{2} {\tilde w_i}^5.\qquad (***1) $$

Clearly, we have:

\bigskip

\underbar {\bf Lemma 3 :}

$$ b_1(y_i,t,\theta)=1-\dfrac{1}{3} Ricci_{y_i}(\theta,\theta)e^{2t}+\ldots, $$

$$ R_g(e^t\theta)=R_g(y_i) + <\nabla R_g(y_i)|\theta > e^t+\dots . $$

According to proposition 1 and lemma 3,

\bigskip

\underbar {\bf Propostion 4 :}

$$ \bar Z_i(\tilde w_i^{\xi_i}-\tilde w_i) \leq V_i {b_1}^{(-2}[(\tilde w_i^{\xi_i})^{5}- \tilde w_i^{5}]+2(\tilde w_i^{\xi_i})^5|V_i^{\xi_i}-V_i|+  $$

$$ +C|e^{2t}-e^{2t^{\xi_i}}|\left [|\nabla_{\theta} {\tilde w_i}^{\xi_i}| + |\nabla_{\theta}^2(\tilde w_i^{\xi_i})|+ |Ricci_{y_i}|[\tilde w_i^{\xi_i}+(\tilde w_i^{\xi_i})^{5}] + |R_g(y_i)| \tilde w_i^{\xi_i} \right ] + C'|e^{3t^{\xi_i}}-e^{3t}|. $$

\underbar {\bf Proof:}

\bigskip

We use proposition 1, we have:

$$ a(y_i,t,\theta)=\log J(y_i,e^t,\theta)=\log b_1, |\partial_t b_1(t)|+|\partial_{tt} b_1(t)|+|\partial_{tt} a(t)|\leq C e^{2t}, $$

and,

$$ |\partial_{\theta_j} b_1|+|\partial_{\theta_j,\theta_k} b_1|+\partial_{t,\theta_j}b_1|+|\partial_{t,\theta_j,\theta_k} b_1|\leq C e^{2t}, $$

then,

$$ |\partial_t b_1(t^{\xi_i})-\partial_t b_1(t)|\leq C'|e^{2t}-e^{2t^{\xi_i}}|,\,\, {\rm sur} \,\, ]-\infty, \log \epsilon_1]\times {\mathbb S}_{2},\forall \,\, x\in B(x_0,\epsilon_1) $$

Locally,

$$ \Delta_{g_{y_i, e^t, {}_{{\mathbb S}_{2}}}}= L_{\theta}(y_i,e^t)=-\dfrac{1}{\sqrt {|{\tilde g}^k(e^t,\theta)|}}\partial_{\theta^l}[{\tilde g}^{\theta^l \theta^j}(e^t,\theta)\sqrt { |{\tilde g}^k(e^t,\theta)|}\partial_{\theta^j}] . $$

Thus, in $ [0,\epsilon_1]\times U^k $, we have,

$$ A_i=\left [{ \left [ \dfrac{1}{\sqrt {|{\tilde g}^k|}}\partial_{\theta^l}({\tilde g}^{\theta^l \theta^j}\sqrt { |{\tilde g}^k|}\partial_{\theta^j}) \right ] }^{\xi_i}- \dfrac{1}{\sqrt {|{\tilde g}^k|}}\partial_{\theta^l}({\tilde g}^{\theta^l \theta^j}\sqrt { |{\tilde g}^k|}\partial_{\theta^j}) \right ](\tilde w_i^{\xi_i}) $$

then, $ A_i=B_i+D_i $ with,

$$ B_i=\left [ {\tilde g}^{\theta^l \theta^j}(e^{t^{\xi_i}}, \theta)-{\tilde g}^{\theta^l \theta^j}(e^t,\theta) \right ] \partial_{\theta^l \theta^j}\tilde w_i^{\xi_i}(t,\theta), $$

and,

$$ D_i=\left [ \dfrac{1}{ \sqrt {| {\tilde g}^k|}(e^{t^{\xi_i}},\theta )  }           \partial_{\theta^l}[{\tilde g }^{\theta^l \theta^j}(e^{t^{\xi_i}},\theta)\sqrt {| {\tilde g}^k|}(e^{t^{\xi_i}},\theta)  ] -\dfrac{1}{ \sqrt {| {\tilde g}^k|}(e^t,\theta) } \partial_{\theta^l} [{\tilde g }^{\theta^l \theta^j}(e^t,\theta)\sqrt {| {\tilde g}^k|}(e^t,\theta) ] \right ] \partial_{\theta^j} \tilde w_i^{\xi_i}(t,\theta), $$

we deduce,

$$ A_i \leq C_k|e^{2t}-e^{2t^{\xi_i}}|\left [ |\nabla_{\theta} \tilde w_i^{\xi_i}| + |\nabla_{\theta}^2(\tilde w_i^{\xi_i})| \right ], $$

If we take $ C=\max \{ C_i, 1 \leq i\leq q \} $ and if w use $ (***1) $, we obtain proposition 4.

We have,

$$ c(y_i,t,\theta)=\left ( \dfrac{n-2}{2} \right )^2+ \dfrac{n-2}{2} \partial_t a + R_g e^{2t}, \qquad (\alpha_1) $$ 

$$ b_2(t,\theta)=\partial_{tt} (\sqrt {b_1})=\dfrac{1}{2 \sqrt { b_1}}\partial_{tt}b_1-\dfrac{1}{4(b_1)^{3/2}}(\partial_t b_1)^2 ,\qquad (\alpha_2) $$ 

$$ c_2=[\dfrac{1}{\sqrt {b_1}} \Delta_{g_{y_i, e^t, {}_{{\mathbb S}_{n-1}}}}(\sqrt{b_1}) + |\nabla_{\theta} \log (\sqrt {b_1})|^2], \qquad (\alpha_3) $$

Then,

$$ \partial_{t}c(y_i,t,\theta)=\dfrac{1}{2}\partial_{tt}a+2e^{2t}R_g(e^t\theta)+e^{3t}<\nabla R_g(e^t\theta)|\theta >, $$

by proposition 1,

$$ |\partial_tc_2|+|\partial_t b_1|+|\partial_t b_2|+|\partial_t c|\leq K_1e^{2t}, $$

Now, we consider the function, $ \bar w_i(t,\theta)=\tilde w_i(t,\theta)-\dfrac{ [u_i(y_i)]^{1/3} \times \min_M u_i}{2} e^{t} $, and $ \lambda >2 >0 $.

\bigskip

For $ t \leq t_i=-(2/3) \log u_i(y_i) $, we have: 

$$ \bar w_i(t,\theta) = e^{t}\left [b_1(t,\theta) e^{-t/2}u_i o \exp_{y_i}(e^t\theta)-\dfrac{ [u_i(y_i)]^{1/3} \times  \min_M u_i}{2} \right ] \geq $$

$$ \geq e^{t} \dfrac{ [u_i(y_i)]^{1/3} \times \min_M u_i}{2}>0, $$

We set, $ \mu_i=\dfrac{[u_i(y_i)]^{1/3} \times \min_M u_i}{2} $.

\bigskip

We use proposition 3 and the same arguments than in [B], we obtain: 

\bigskip

\underbar {\bf Lemma 4:}

\smallskip

There exists $ \nu <0 $ such that for $ \mu \leq \nu $ :

$$ \bar w_i^{\mu}(t,\theta)-\bar w_i(t,\theta) \leq 0, \,\, \forall \,\, (t,\theta) \in [\mu,t_i] \times {\mathbb S}_{2}, $$

We set, $ \lambda_i=-2 \log u_i(y_i) $, then, 

\bigskip

\underbar {\bf Lemma 5:}

$$ \bar w_i(\lambda_i,\theta)-\bar w_i(\lambda_i+4,\theta) >0. $$

\underbar {\bf Proof of lemma 5:}

\bigskip

Clearly:

$$ \bar w_i(\lambda_i,\theta)-\bar w_i(\lambda_i+4,\theta)= \tilde w_i(\lambda_i,\theta)-\tilde w_i(\lambda_i+4,\theta)+\mu_i(e^4-1), $$

we deduce lemma 6 from proposition 3.

\bigskip

Let, $ \xi_i=\sup \{ \mu \leq \lambda_i+2, \bar w_i^{\xi_i}(t,\theta)-\bar w_i(t,\theta) \leq 0, \,\, \forall \,\, (t,\theta)\in [\xi_i,t_i]\times {\mathbb S}_{2} \} $.

\bigskip

The real $ \xi_i $ exists (see [B]), if we use $ (***2) $, we have:

$$ \tilde w_i^{\xi_i}(t,\theta) + |\nabla_{\theta} \tilde  w_i^{\xi_i}(t,\theta)|+|\nabla_{\theta}^2  \tilde w_i^{\xi_i}(t,\theta)|\leq C(R),\,\,\, \forall \,\, (t,\theta) \in ]-\infty,\log R]\times {\mathbb S}_{2}, $$

We can write:

$$ \bar Z_i(\bar w_i^{\xi_i}-\bar w_i)=\bar Z_i(\tilde w_i^{\xi_i}-\tilde w_i)-\mu_i \bar Z_i(e^{2t^{\xi_i}}-e^{2t}), $$

$$ -\bar Z_i(e^{2t^{\xi_i}}-e^{2t})=[1-\dfrac{1}{4}- \dfrac{3}{2}\partial_t a-R_ge^{2t}+b_1^{-1/2} b_2-c_2](e^{2t^{\xi_i}}-e^{2t}) \leq c_1(e^{t^{\xi_i}}-e^{t}), $$

with $ c_1>0 $, car $ |\partial_t a|+|\partial_t b_1|+|\partial_{tt} b_1|+|\partial_{t,\theta_j} b_1|+|\partial_{t,\theta_j,\theta_k} b_1|\leq C'e^{2t}<1 $, for $ t $ very small.

\bigskip

We use proposition 4, to obtain on, $ [\xi_i,t_i]\times {\mathbb S}_{2} $,

$$ \bar Z_i(\bar w_i^{\xi_i}-\bar w_i) \leq c_2[V_i(\tilde w_i^{\xi_i})^{5}-\tilde w_i^{5}]+|V_i^{\xi_i}-V_i|(w_i^{\xi_i})^5 + [\mu_i c_1- C'(R)](e^{t^{\xi_i}}-e^{t}) \leq 0, $$

\bigskip

Like in [B], after using Hopf maximum principle, we have,

$$ \sup_{\theta \in {\mathbb S}_{2}} [\bar w_i^{\xi_i}(t_i,\theta)-\bar w_i(t_i,\theta)]=0. $$

We have:

$$ \bar w_i^{\xi_i}(t_i,\theta_i)-\bar w_i(t_i,\theta_i)=0, \,\, \forall \,\, i. $$

We deduce,

$$ [u_(y_i)]^{1/3} \times \inf_M u_i \leq c,\,\, \forall \,\, i. $$

\bigskip

\underbar {\bf References:}

\bigskip

[Au] T. Aubin. Some Nonlinear Problems in Riemannian Geometry. Springer-Verlag 1998.

\smallskip

[B] S.S Bahoura. Majorations du type $ \sup u \times \inf u \leq c $ pour l'\'equation de la courbure scalaire sur un ouvert de $ {\mathbb R}^n, n\geq 3 $. J. Math. Pures. Appl.(9) 83 2004 no, 9, 1109-1150.

\smallskip

[B-L-S] H. Brezis, Yy. Li Y-Y, I. Shafrir. A sup+inf inequality for some
nonlinear elliptic equations involving exponential
nonlinearities. J.Funct.Anal.115 (1993) 344-358.

\smallskip

[BM] H.Brezis and F.Merle, Uniform estimates and blow-up bihavior for solutions of $ -\Delta u=Ve^u $ in two dimensions, Commun Partial Differential Equations 16 (1991), 1223-1253.

\smallskip

[C-G-S], L. Caffarelli, B. Gidas, J. Spruck. Asymptotic symmetry and local
behavior of semilinear elliptic equations with critical Sobolev
growth. Comm. Pure Appl. Math. 37 (1984) 369-402.

\smallskip

[CL1] C-C.Chen, C-S. Lin. Estimates of the conformal scalar curvature
equation via the method of moving planes. Comm. Pure
Appl. Math. L(1997) 0971-1017.

\smallskip

[CL2] A sharp sup+inf inequality for a nonlinear elliptic equation in ${\mathbb R}^2$.
Commun. Anal. Geom. 6, No.1, 1-19 (1998).

\smallskip

[L,P] J.M. Lee, T.H. Parker. The Yamabe problem. Bull.Amer.Math.Soc (N.S) 17 (1987), no.1, 37 -91.

\smallskip

[L1] Yy. Li. Prescribing scalar curvature on $ {\mathbb S}_n $ and related
Problems. C.R. Acad. Sci. Paris 317 (1993) 159-164. Part
I: J. Differ. Equations 120 (1995) 319-410. Part II: Existence and
compactness. Comm. Pure Appl.Math.49 (1996) 541-597.

\smallskip

[L2] Yy. Li. Harnack Type Inequality: the Method of Moving Planes. Commun. Math. Phys. 200,421-444 (1999).

\smallskip

[L-Z] Yy. Li, L. Zhang. A Harnack type inequality for the Yamabe equation in low dimensions.  Calc. Var. Partial Differential Equations  20  (2004),  no. 2, 133--151.

\smallskip

[L-Zh] Yy.Li, M. Zhu. Yamabe Type Equations On Three Dimensional Riemannian Manifolds. Commun.Contem.Mathematics, vol 1. No.1 (1999) 1-50.

\smallskip

[S] I. Shafrir. A sup+inf inequality for the equation $ -\Delta u=Ve^u $. C. R. Acad.Sci. Paris S\'er. I Math. 315 (1992), no. 2, 159-164.

\end{document}